\pdfoutput=1
\RequirePackage{ifpdf}
\ifpdf 
\documentclass[pdftex]{sigma}
\else
\documentclass{sigma}
\fi

\newcommand{\bn}{{\mathbb{N}}}
\newcommand{\br}{{\mathbb{R}}}

\newcommand{\bz}{{\mathbb{Z}}}
\newcommand{\bc}{{\mathbb{C}}}

\newcommand{\bt}{{\mathbb{T}}}

\newcommand{\ca}{{\mathcal{A}}}
\newcommand{\cb}{{\mathcal{B}}}
\newcommand{\cd}{{\mathcal{D}}}
\newcommand{\cf}{{\mathcal{F}}}
\newcommand{\ch}{{\mathcal{H}}}
\newcommand{\cg}{{\mathcal{G}}}
\newcommand{\cl}{{\mathcal{L}}}

\newcommand{\css}{{\mathcal{S}}}

\renewcommand{\a}{\alpha}
\renewcommand{\b}{\beta}
\renewcommand{\l}{\lambda}
\renewcommand{\ll}{\Lambda}
\newcommand{\s}{\sigma}

\newcommand{\ep}{\varepsilon}
\newcommand{\z}{\zeta}

\newcommand{\ovl}{\overline}

\newcommand{\wt}{\widetilde}

\DeclareMathOperator{\diag}{\rm diag}

\DeclareMathOperator{\tr}{tr}

\numberwithin{equation}{section}

\newtheorem{Theorem}{Theorem}[section]
\newtheorem{Corollary}[Theorem]{Corollary}
\theoremstyle{definition}
\newtheorem{Remark}[Theorem]{Remark}
\newtheorem{Example}[Theorem]{Example}

\begin{document}

\newcommand{\arXivNumber}{1809.07136}

\renewcommand{\PaperNumber}{050}

\FirstPageHeading
	
\ShortArticleName{On Direct Integral Expansion for Periodic Block-Operator Jacobi Matrices}
	
\ArticleName{On Direct Integral Expansion for Periodic\\ Block-Operator Jacobi Matrices and Applications}
	
\Author{Leonid GOLINSKII~$^\dag$ and Anton KUTSENKO~$^\ddag$}
	
\AuthorNameForHeading{L.~Golinskii and A.~Kutsenko}
	
\Address{$^\dag$~B.~Verkin Institute for Low Temperature Physics and Engineering,\\
\hphantom{$^\dag$}~47 Science Ave., Kharkiv 61103, Ukraine}
\EmailD{\href{mailto:golinskii@ilt.kharkov.ua}{golinskii@ilt.kharkov.ua}}
	
\Address{$^\ddag$~Jacobs University, Campus Ring 1, 28759 Bremen, Germany}
\EmailD{\href{mailto:akucenko@gmail.com}{akucenko@gmail.com}}

\ArticleDates{Received December 03, 2018, in final form June 23, 2019; Published online July 02, 2019}
	
\Abstract{We construct a functional model (direct integral expansion) and study the spectra of certain periodic block-operator Jacobi matrices, in particular, of general 2D partial difference operators of the second order. We obtain the upper bound, optimal in a sense, for the Lebesgue measure of their spectra. The examples of the operators for which there are several gaps in the spectrum are given.}
	
	
\Classification{47B36; 47B39; 35P15}

\section{Introduction}

Consider a block-operator Jacobi matrix on the Hilbert space $\cf=\ell^2(\bz, \ch)$
\begin{gather}\label{genbloper}
(Ju)_q:=A_{q-1}u_{q-1}+B_qu_q+A_q^*u_{q+1}, \qquad B_q=B_q^*, \qquad q\in\bz, \qquad u_q\in\ch,
\end{gather}
where the blocks $\{A_q, B_q\}$ are bounded linear operators on the Hilbert space~$\ch$. Under the standard assumption $\sup\limits_{q\in\bz}(\|A_q\|+\|B_q\|)<\infty$ on the entries, the matrix $J$ generates a bounded
and self-adjoint operator $J$ on $\cf$.

We are primarily interested in the case when $\ch=\ell^2(\bz)$, and the blocks $B_q$, $A_q$ are themselves Hermitian (not necessarily real symmetric) 1D Jacobi matrices
\begin{gather}\label{bloper}
(Ju)_q:=A_{q-1}u_{q-1}+B_qu_q+A_qu_{q+1}, \qquad q\in\bz, \qquad u_q=(u_q(l))_l\in\ell^2(\bz).
\end{gather}
The block matrices of this type are known as Jacobi-block-Jacobi matrices. So,
\begin{gather}\label{jacbloc}
A_q=J\bigl(\{a_{q,r}\}_{r\in\bz}, \{\a_{q,r}\}_{r\in\bz}\bigr), \qquad B_q=J\bigl(\{c_{q,r}\}_{r\in\bz}, \{b_{q,r}\}_{r\in\bz}\bigr),
\end{gather}
with real entries $\{a_{q,r}\}$ $(\{c_{q,r}\})$ along the main diagonal, and complex entries $\{\a_{q,r}\}$ $(\{b_{q,r}\})$ and their conjugates along the diagonals above and below the main one, respectively.

Such block-operator model arises when we deal with 2D partial difference operators $L$ of the second order of the form
\begin{gather}\label{partdif}
(L\wt u)_{ik} :=a_{i-1,k}u_{i-1,k}+a_{i,k}u_{i+1,k}+b_{i,k-1}u_{i,k-1}+b_{i,k}u_{i,k+1}+c_{i,k}u_{i,k},
\end{gather}
on the Hilbert space $\ell^2\big(\bz^2\big)$. A natural isometry $U_0$ between $\ell^2\big(\bz^2\big)$ and $\cf$
\begin{gather*}
\wt u=(u_{i,k})_{i,k\in\bz} \rightarrow u=U_0\wt u=(u_j)_{j\in\bz}\colon \ u_j=(\dots, u_{j,-1}, u_{j,0}, u_{j,1},\dots)'
\end{gather*}
transforms the operator $L$ in \eqref{partdif} into $J:=U_0 LU_0^{-1}$ on $\cf$, where
\begin{gather}
(J u)_j = A_{j-1}u_{j-1}+B_j u_j+A_j u_{j+1}, \qquad j\in\bz, \nonumber\\
A_j :=\diag (a_{j,k})_{k\in\bz}, \qquad B_j:=J\bigl(\{c_{j,k}\}_{k\in\bz}, \{b_{j,k}\}_{k\in\bz}\bigr),\label{opmod}
\end{gather}
see, e.g., \cite[Section VII.3]{B68}. This is a special case of the model operator $J$ \eqref{bloper} with $\a_{j,k}=0$.
In particular, for 2D discrete Schr\"odinger operators $H$ ($a_{j,k}=b_{j,k}=1$) we have
\begin{gather}\label{disshr}
(H\wt u)_{i,k} :=u_{i-1,k}+u_{i+1,k}+u_{i,k-1}+u_{i,k+1}+c_{i,k}u_{i,k},
\end{gather}
so
\begin{gather*} A_j=I, \qquad B_j=J\bigl(\{c_{j,k}\}_{k\in\bz}, \{1\}\bigr), \end{gather*}
1D discrete Schr\"odinger operators.

Let $p_j\in\bn$, $j=1,2$. A sequence of complex numbers $\{w_{q,r}\}_{q,r\in\bz}$ is called $(p_1,p_2)$-periodic if
\begin{gather*}
w_{q+k_1p_1, r+k_2p_2}=w_{q,r}, \qquad \forall q,r,k_1,k_2\in\bz.
\end{gather*}
The operator $J$ in \eqref{bloper} is called $(p_1,p_2)$-periodic if $\{a_{q,r}\}$, $\{\a_{q,r}\}$, $\{c_{q,r}\}$, and $\{b_{q,r}\}$, $q,r\in\bz$, are $(p_1,p_2)$-periodic. Equivalently,
\begin{enumerate}\itemsep=0pt
 \item[1)] all the blocks $A_q$, $B_q$ are $p_2$-periodic 1D Jacobi matrices;
 \item[2)] $J$ is block periodic with period $p_1$: $A_{q+p_1}=A_q$, $B_{q+p_1}=B_q$.
\end{enumerate}
We say that the partial difference operator $L$ \eqref{partdif} is $(p_1,p_2)$-periodic, if all the coefficients are $(p_1,p_2)$-periodic,
or equivalently, $J$ \eqref{opmod} is $(p_1,p_2)$-periodic.

In Section \ref{s1} we recall a direct integral expansion (a functional model) for the $(p_1,p_2)$-periodic operator
$J$ in \eqref{bloper}--\eqref{jacbloc} and establish the banded structure of its spectrum $\s(J)$. In Section~\ref{s2} we estimate the
Lebesgue measure of $\s(J)$.

\begin{Theorem}\label{upbound}
Let $p_1,p_2\ge3$. The Lebesgue measure of the spectrum for the periodic ope\-ra\-tor~$J$ in \eqref{bloper}--\eqref{jacbloc} admits the upper bound
\begin{gather}\label{upboun}
|\s(J)|\le \min_{(m,n)\in\bz^2}R_{m,n},
\end{gather}
where
\begin{gather}\label{upperbound}
R_{m,n}:=4\sum_{j=1}^{p_1} (|b_{j,n}|+2|\a_{j,n}|)+4\sum_{k=1}^{p_2} (|a_{m,k}|+2|\a_{m,k}|)-8|\a_{m,n}|.
\end{gather}
\end{Theorem}

The sequence $\{R_{m,n}\}$ is $(p_1,p_2)$-periodic, so minimum in \eqref{upboun} is actually taken over the finite set of indices
$m=1,\dots,p_1$ and $n=1,\dots,p_2$.

The similar result holds for $\min(p_1,p_2)=2$, see Remark~\ref{case2} below.

Note that there is a simple general bound for $|\s(J)|$ (which has nothing to do with periodicity) based on the fact that $J$ is
a three-diagonal block-matrix
\begin{gather}
|\s(J)| \le 2\|J\|\le 2(2\max_n\|A_n\| +\max_n\|B_n\|) \nonumber\\
\hphantom{|\s(J)|}{} \le \max_n (8\max_m|\a_{n,m}|+4\max_m|a_{n,m}|) +\max_n(4\max_m |b_{n,m}|+2\max_m|c_{n,m}|).\label{upboun1}
\end{gather}
The point is that certain parameters, such as $c_{n,m}$, which appear in \eqref{upboun1}, do not enter~\eqref{upboun}. So, once some values of $c_{n,m}$
are large enough, bound \eqref{upboun} is better that~\eqref{upboun1}. On the other hand, \eqref{upboun} contains sums of the
entries compared to~\eqref{upboun1}, which does not.

There is another upper bound for the length of the spectrum, based on Gershgorin's theorem, which is discussed in Remark~\ref{gersh}.

For the 1D scalar, $p$-periodic Jacobi operator $J$ the estimates for the spectrum
\begin{gather*} |\s(J)|\le 4(a_1a_2\cdots a_p)^{1/p} \end{gather*}
were obtained in \cite{KoKr, AK06}. Here $a_n$ are the off-diagonal entries of $J$. Recently the second author~\cite{AK10} improved this result to
\begin{gather*} |\s(J)|\le 4\min_n a_n. \end{gather*}
We see that both bounds do not depend on the diagonal entries, so the lack of $c_{m,n}$ in~\eqref{upboun} looks reasonable.

\begin{Corollary}\label{cor02}
For a $2D$ periodic, partial difference operators of the second order $L$ in~\eqref{partdif} the spectral estimate is
\begin{gather*} \frac14 |\s(L)|\le \min_n \sum_{j=1}^{p_1} |b_{j,n}|+\min_m \sum_{k=1}^{p_2} |a_{m,k}|. \end{gather*}
In particular, $|\s(H)|\le4(p_1+p_2)$ for $2D$ periodic, discrete Schr\"odinger operators~$H$~\eqref{disshr}.
\end{Corollary}

The estimate $|\s(H)|\le4(p_1+p_2)$ was previously obtained in \cite{AK15}.

We complete the paper with examples of 2D discrete Schr\"odinger operators with explicitly computed spectra, and an operator, which has a maximal number of gaps in its spectrum.

\section{Direct integral expansion}\label{s1}

We begin with auxiliary, Hermitian matrix-functions of the order $p_2$
\begin{gather}\label{aux1}
\ca_n(x_2)=\begin{bmatrix}
a_{n,1} & \a_{n,1} & & & & {\rm e}^{{\rm i}x_2}\ovl{\a}_{n,p_2} \\
\ovl{\a}_{n,1} & a_{n,2} & \a_{n,2} & & & \\
& \ovl{\a}_{n,2} & a_{n,3} & \a_{n,3} & & \\
& & \ddots & \ddots & \ddots & \\
& & & \ovl{\a}_{n,p_2-2} & a_{n,p_2-1} & \a_{n,p_2-1} \\
{\rm e}^{-{\rm i}x_2}\a_{n,p_2} & & & & \ovl{\a}_{n,p_2-1} & a_{n,p_2}
\end{bmatrix}
\end{gather}
and
\begin{gather}\label{aux2}
\cb_n(x_2)=\begin{bmatrix}
c_{n,1} & b_{n,1} & & & & {\rm e}^{{\rm i}x_2}\ovl{b}_{n,p_2} \\
\ovl{b}_{n,1} & c_{n,2} & b_{n,2} & & & \\
& \ovl{b}_{n,2} & c_{n,3} & b_{n,3} & & \\
& & \ddots & \ddots & \ddots & \\
& & & \ovl{b}_{n,p_2-2} & c_{n,p_2-1} & b_{n,p_2-1} \\
{\rm e}^{-{\rm i}x_2}b_{n,p_2} & & & & \ovl{b}_{n,p_2-1} & c_{n,p_2}
\end{bmatrix},
\end{gather}
$x_2\in[0,2\pi)$, and combine them in a single block matrix-function $\css$ of the order $p:=p_1p_2$,{\samepage
\begin{gather}\label{symb}
\css=\css(x_1,x_2)=\begin{bmatrix}
\cb_1 & \ca_1 & & & & {\rm e}^{{\rm i}x_1}\ca_{p_1} \\
\ca_1 & \cb_2 & \ca_2 & & & \\
& \ca_2 & \cb_3 & \ca_3 & & \\
& & \ddots & \ddots & \ddots & \\
& & & \ca_{p_1-2} & \cb_{p_1-1} & \ca_{p_1-1} \\
{\rm e}^{-{\rm i}x_1}\ca_{p_1} & & & & \ca_{p_1-1} & \cb_{p_1}
\end{bmatrix},
\end{gather}
$x_1,x_2\in[0,2\pi)$, the main object under consideration, known as a {\it symbol}.}

Denote $\mathbb{T}^2:=[0,2\pi)\times [0,2\pi)$, and put
\begin{gather*}
\cl = \int_{\mathbb{T}^2}^{\oplus} \bc^p \frac{{\rm d}x_1{\rm d}x_2}{4\pi^2}=L^2\big(\mathbb{T}^2, \bc^p\big)
= \left\{
\begin{bmatrix}
g_1(x_1,x_2) \\
g_2(x_1,x_2) \\
\vdots \\
g_{p_1}(x_1,x_2)
\end{bmatrix}\colon g_k(x_1,x_2)=
\begin{bmatrix}
g_{k,1}(x_1,x_2) \\
g_{k,2}(x_1,x_2) \\
\vdots \\
g_{k,p_2}(x_1,x_2)
\end{bmatrix} \right\},
\end{gather*}
where
\begin{gather*}
g_{k,j}(x_1,x_2) =\sum_{m,n\in\bz}\hat g_{k,j}(m,n){\rm e}^{{\rm i}mx_1+{\rm i}nx_2}\in L^2\big(\mathbb{T}^2\big), \nonumber\\
k =1,2,\dots, p_1, \qquad j=1,2,\dots,p_2.
\end{gather*}
We have
\begin{gather*}
\|g\|^2_\cl=\sum_{k,j}\|g_{k,j}\|^2_{L^2(\mathbb{T}^2)}=\sum_{k,j,m,n}|\hat g_{k,j}(m,n)|^2.
\end{gather*}

There is a natural isometry $U\colon \cl\to\cf=\ell^2\big(\bz, \ell^2(\bz)\big)$ which acts by
\begin{gather}
U
\begin{bmatrix}
g_1 \\
g_2 \\
\vdots \\
g_{p_1}
\end{bmatrix}=u=(u_r)_{r\in\bz}, \qquad u_r =(u_r(s))_{s\in\bz}\colon \ u_{k+p_1 m}(j+p_2 n) =\hat g_{k,j}(m,n).\label{unitary}
\end{gather}
Assume that $p_1, p_2\ge3$. The above symbol $\css$ defines a multiplication operator $M(\css)$ on $\cl$ by
\begin{gather}\label{mult1}
M(\css)g=\css
\begin{bmatrix}
g_1(x_1,x_2) \\
g_2(x_1,x_2) \\
\vdots \\
g_{p_1}(x_1,x_2)
\end{bmatrix}=
\begin{bmatrix}
h_1(x_1,x_2) \\
h_2(x_1,x_2) \\
\vdots \\
h_{p_1}(x_1,x_2)
\end{bmatrix}
\end{gather}
with
\begin{gather}
h_1(x_1,x_2) =\cb_1 g_1+\ca_1 g_2+{\rm e}^{{\rm i}x_1}\ca_{p_1}g_{p_1},\nonumber\\
h_l(x_1,x_2) = \ca_{l-1}g_{l-1}+\cb_l g_l+\ca_l g_{l+1}, \qquad l=2,\dots,p_1-1, \nonumber\\
h_{p_1}(x_1,x_2) = \ca_{p_1-1}g_{p_1-1}+\cb_{p_1} g_{p_1}+{\rm e}^{-{\rm i}x_1}\ca_{p_1} g_1.\label{mult2}
\end{gather}

The result below is a cornerstone of 2D discrete, Floquet--Bloch theory, see \cite[Section~5.3]{SSz} for 1D theory (the calculation in dimension $2$ is completely analogous).

\begin{Theorem}Let $p_1,p_2\ge3$. The $(p_1,p_2)$-periodic operator $J$ \eqref{bloper}--\eqref{jacbloc} is unitarily equivalent to the multiplication
operator $M(\css)$ \eqref{mult1}--\eqref{mult2}
\begin{gather*} J=UM(\css)U^{-1}, \end{gather*}
$U$ is defined in \eqref{unitary}.
\end{Theorem}

\begin{Remark}\label{smallperiod}
There is nothing special in the case $\min(p_1,p_2)=2$, but the symbol looks differently. Precisely, if $p_1\ge3$, $p_2=2$, we have
\begin{gather}
\ca_n =
\begin{bmatrix}
a_{n,1} & \a_{n,1}+{\rm e}^{{\rm i}x_2}\ovl{\a}_{n,2} \\
\ovl{\a}_{n,1}+{\rm e}^{-{\rm i}x_2}\a_{n,2} & a_{n,2}
\end{bmatrix}, \nonumber \\
\cb_n =
\begin{bmatrix}
c_{n,1} & b_{n,1}+{\rm e}^{{\rm i}x_2}\ovl{b}_{n,2} \\
\ovl{b}_{n,1}+{\rm e}^{-{\rm i}x_2}b_{n,2} & c_{n,2}
\end{bmatrix}, \qquad n=1,\dots,p_1,\label{smalper}
\end{gather}
and $\css$ is of the form \eqref{symb}. If $p_1=2$, $p_2\ge3$, we take $\ca_1$, $\cb_1$, $\ca_2$, $\cb_2$ as in~\eqref{aux1}, \eqref{aux2}, and
\begin{gather}\label{smalpersym}
\css (x_1,x_2)=
\begin{bmatrix}
\cb_1 & \ca_{1}+{\rm e}^{{\rm i}x_1}\ca_{2} \\
\ca_{1}+{\rm e}^{-{\rm i}x_1}\ca_{2} & \cb_2
\end{bmatrix}.
\end{gather}
Finally, if $p_1=p_2=2$, then $\ca_n$, $\cb_n$, and $\css$ are of the form \eqref{smalper}--\eqref{smalpersym}.
\end{Remark}

Denote by $\l_1\ge\l_2\ge\dots\ge\l_p$, $\l_j=\l_j(x_1,x_2)$, the set of all eigenvalues of $\css$, labeled in the non-increasing order. According to the general result on the spectrum of multiplication operators,
\begin{gather}
\s(J) =\s(M(\css))=\bigcup_{k=1}^p \ll_k, \qquad
\ll_k :=\big[\inf_{\mathbb{T}^2}\l_k(x_1,x_2), \sup_{\mathbb{T}^2}\l_k(x_1,x_2)\big],\label{spect}
\end{gather}
the $k$'s band in the spectrum. So we come to the following

\begin{Corollary}\label{spectr}The spectrum of the $(p_1,p_2)$-periodic Jacobi-block-Jacobi matrix \eqref{bloper} has the banded structure
\begin{gather*}
\s(J)=\bigcup_{k=1}^p \ll_k, \qquad p:=p_1p_2,
\end{gather*}
with the closed intervals $\ll_k$~\eqref{spect}. So, the number of the gaps in the spectrum does not exceed $p-1$.
\end{Corollary}

Note that for a $(p_1,p_2)$-periodic, 2D discrete Schr\"odinger operator the symbol $\css$ takes the form
\begin{gather*}
\css(x_1,x_2) =\begin{bmatrix}
\cb_1(x_2) & I_{p_2} & & & & {\rm e}^{{\rm i}x_1}I_{p_2} \\
I_{p_2} & \cb_2(x_2) & I_{p_2} & & & \\
& I_{p_2} & \cb_3(x_2) & I_{p_2} & & \\
& & \ddots & \ddots & \ddots & \\
& & & I_{p_2} & \cb_{p_1-1}(x_2) & I_{p_2} \\
{\rm e}^{-{\rm i}x_1}I_{p_2} & & & & I_{p_2} & \cb_{p_1}(x_2)
\end{bmatrix}, \nonumber\\
\cb_n(x_2) =\begin{bmatrix}
c_{n,1} & 1 & & & & {\rm e}^{{\rm i}x_2} \\
1 & c_{n,2} & 1 & & & \\
& 1 & c_{n,3} & 1 & & \\
& & \ddots & \ddots & \ddots & \\
& & & 1 & c_{n,p_2-1} & 1 \\
{\rm e}^{-{\rm i}x_2} & & & & 1 & c_{n,p_2}
\end{bmatrix}, \qquad x_1,x_2\in[0,2\pi).
\end{gather*}

\section{Spectral estimates for periodic block-Jacobi operators}\label{s2}

In this section we are aimed at proving Theorem \ref{upbound}. By Corollary \ref{spectr},
\begin{gather}\label{spband}
\s(J)=\s(M(\css))=\bigcup_{k=1}^p \ll_k, \qquad \ll_k=[l_k,r_k]
\end{gather}
are the closed intervals, swept by the $k$-th eigenvalue $\l_k(x_1,x_2)$, arranged in the non-increasing order, as the pair
$(x_1,x_2)$ runs over $\mathbb{T}^2$.

We are looking for two constant matrices $\css_{\pm}$, i.e., independent of $(x_1,x_2)$, so that
\begin{gather}\label{eigenboun}
\css_-\le\css\le\css_+ \quad \Rightarrow \quad \l_k^-\le\l_k(x_1,x_2)\le\l_k^+,
\end{gather}
where $\l_1^{\pm}\ge\l_2^{\pm}\ge\dots\ge\l_p^{\pm}$ are eigenvalues of $\css_\pm$, arranged in the non-increasing order. Hence,
\begin{gather}\label{spectrboun1}
|\s(M(\css))|\le\sum_{k=1}^p (r_k-l_k)\le\sum_{k=1}^p (\l_k^+-\l_k^-)=\tr(\css_+-\css_-).
\end{gather}
To this end put
\begin{gather*}
\ca_n^1 =\begin{bmatrix}
a_{n,1} & \a_{n,1} & & & & \\
\ovl{\a}_{n,1} & a_{n,2} & \a_{n,2} & & & \\
& \ovl{\a}_{n,2} & a_{n,3} & \a_{n,3} & & \\
& & \ddots & \ddots & \ddots & \\
& & & \ovl{\a}_{n,p_2-2} & a_{n,p_2-1} & \a_{n,p_2-1} \\
& & & & \ovl{\a}_{n,p_2-1} & a_{n,p_2}
\end{bmatrix}, \\
\ca_n^2(x_2)=\begin{bmatrix}
& & & & & {\rm e}^{{\rm i}x_2}\ovl{\a}_{n,p_2} \\
& & & & & \\
& & & & & \\
& & & & & \\
& & & & & \\
{\rm e}^{-{\rm i}x_2}\a_{n,p_2} & & & & &
\end{bmatrix},
\end{gather*}
so $\ca_n^1$ and $\ca_n^2$ are Hermitian matrices, and let $\ca_n=\ca_n^1+\ca_n^2$. Similarly,
\begin{gather*}
\cb_n^1 =\begin{bmatrix}
c_{n,1} & b_{n,1} & & & & \\
\ovl{b}_{n,1} & c_{n,2} & b_{n,2} & & & \\
& \ovl{b}_{n,2} & c_{n,3} & b_{n,3} & & \\
& & \ddots & \ddots & \ddots & \\
& & & \ovl{b}_{n,p_2-2} & c_{n,p_2-1} & b_{n,p_2-1} \\
& & & & \ovl{b}_{n,p_2-1} & c_{n,p_2}
\end{bmatrix}, \\
\cb_n^2(x_2) =\begin{bmatrix}
 & & & & & {\rm e}^{{\rm i}x_2}\ovl{b}_{n,p_2} \\
 & & & & & \\
& & & & & \\
& & & & & \\
& & & & & \\
{\rm e}^{-{\rm i}x_2}b_{n,p_2} & & & & &
\end{bmatrix},
\end{gather*}
so $\cb_n^1$ and $\cb_n^2$ are Hermitian matrices, and let $\cb_n=\cb_n^1+\cb_n^2$. It follows that
\begin{gather}\label{spli}
\css(x_1,x_2)=\css_1+\css_2(x_1,x_2)
\end{gather}
with
\begin{gather*}\allowdisplaybreaks
\css_1 =\begin{bmatrix}
\cb_1^1 & \ca_1^1 & & & & \\
\ca_1^1 & \cb_2^1 & \ca_2^1 & & & \\
& \ca_2^1 & \cb_3^1 & \ca_3^1 & & \\
& & \ddots & \ddots & \ddots & \\
& & & \ca_{p_1-2}^1 & \cb_{p_1-1}^1 & \ca_{p_1-1}^1 \\
& & & & \ca_{p_1-1}^1 & \cb_{p_1}^1
\end{bmatrix}, \\
\css_2(x_1,x_2) =\begin{bmatrix}
\cb_1^2 & \ca_1^2 & & & & {\rm e}^{{\rm i}x_1}\ca_{p_1} \\
\ca_1^2 & \cb_2^2 & \ca_2^2 & & & \\
& \ca_2^2 & \cb_3^2 & \ca_3^2 & & \\
& & \ddots & \ddots & \ddots & \\
& & & \ca_{p_1-2}^2 & \cb_{p_1-1}^2 & \ca_{p_1-1}^2 \\
{\rm e}^{-{\rm i}x_1}\ca_{p_1} & & & & \ca_{p_1-1}^2 & \cb_{p_1}^2
\end{bmatrix}.
\end{gather*}
In decomposition \eqref{spli} of the symbol, $\css_1$ is a constant matrix, and $\css_2$ is sparse, i.e., it contains few nonzero entries.

Next,
\begin{gather*}
\css_2=\css_{3}+\css_{4}+\sum_{j=1}^{p_1-1} E_j
\end{gather*}
with
\begin{gather*} \css_{3}:=\diag\big(\cb_1^2, \dots, \cb_{p_1}^2\big), \qquad
\css_{4}:=
\begin{bmatrix}
& & & & {\rm e}^{{\rm i}x_1}\ca_{p_1} \\
& & & & & \\
& & & & & \\
& & & & & \\
{\rm e}^{-{\rm i}x_1}\ca_{p_1} & & & & &
\end{bmatrix},
\end{gather*}
and
\begin{gather*}
E_1 :=
\begin{bmatrix}
0 & \ca_1^2 & 0 & & \\
\ca_1^2 & 0 & 0 & \ddots & & \\
0 & 0 & 0 & \ddots & & \\
& \ddots & \ddots & \ddots & \ddots & \\
& & & & &
\end{bmatrix}, \qquad
E_2:=
\begin{bmatrix}
0 & 0 & 0 & & \\
0 & 0 & \ca_2^2 & \ddots & & \\
0 & \ca_2^2 & 0 & \ddots & & \\
& \ddots & \ddots & \ddots & \ddots & \\
& & & & &
\end{bmatrix}, \qquad \dots \\
E_{p_1-1} :=
\begin{bmatrix}
& & & & \\
& \ddots & \ddots & \ddots & \ddots & \\
& \ddots & 0 & 0 & 0 & \\
& \ddots & 0 & 0 & \ca_{p_1-1}^2 & \\
& \ddots & 0 & \ca_{p_1-1}^2 & 0
\end{bmatrix}.
\end{gather*}
Since $\ca_{p_1}=G_1+G_1^*+G_2$ with
\begin{gather*} G_1:=
\begin{bmatrix}
0 & & & & {\rm e}^{{\rm i}x_2}\ovl{\a}_{p_1,p_2} \\
\ovl{\a}_{p_1,1} & 0 & & & \\
& \ovl{\a}_{p_1,2} & 0 & & \\
& & \ddots & \ddots & \\
& & & \ovl{\a}_{p_1,p_2-1} & 0
\end{bmatrix}, \qquad
G_2:=\diag(a_{p_1,1}, \dots, a_{p_1,p_2}),
\end{gather*}
the matrices of the order $p_2$, we have
\begin{gather*}
\css_{4} =
\begin{bmatrix}
& & & {\rm e}^{{\rm i}x_1}G_1 \\
& & & \\
{\rm e}^{-{\rm i}x_1}G_1^* & & &
\end{bmatrix}+
\begin{bmatrix}
& & {\rm e}^{{\rm i}x_1}G_1^* \\
& & \\
{\rm e}^{-{\rm i}x_1}G_1 & &
\end{bmatrix}+
\begin{bmatrix}
& & {\rm e}^{{\rm i}x_1}G_2 \\
& & \\
{\rm e}^{-{\rm i}x_1}G_2 & &
\end{bmatrix} \\
\hphantom{\css_{4}}{} =\css_{4}'+\css_{4}''+\css_{4}''',
\end{gather*}
so
\begin{gather}\label{decom2}
\css(x_1,x_2)=\css_1+\css_{3}+\css_{4}'+\css_{4}''+\css_{4}'''+\sum_{j=1}^{p_1-1} E_j.
\end{gather}

As the next step toward \eqref{eigenboun}, let us turn to the absolute value $|A|=(A^*A)^{1/2}$ of a matrix $A$, which occurs in the polar representation $A=V|A|$. For an Hermitian matrix $A$, its absolute value can be defined as follows. Denote by $P_+$ ($P_-$) the projection onto the nonnegative (negative) eigenspace of~$A$. Then
\begin{gather*} A=A_+-A_-, \qquad |A|=A_++A_-, \qquad A_{\pm}:=\pm P_{\pm}A\ge0. \end{gather*}
It is clear from this definition that
\begin{gather}\label{absval}
-|A|\le A\le |A|, \qquad -\sum_{j=1}^N |A_j|\le \sum_{j=1}^N A_j\le \sum_{j=1}^N |A_j|.
\end{gather}
We apply \eqref{absval} to decomposition \eqref{decom2} to obtain \eqref{eigenboun} with
\begin{gather*}
\css_\pm:=\css_1\pm D, \qquad D:=|\css_{3}|+|\css_{4}'|+|\css_{4}''|+|\css_{4}'''|+\sum_{j=1}^{p_1-1} |E_j|,
\end{gather*}
and so, by \eqref{spectrboun1},
\begin{gather}\label{specboun}
|\s(J)|=|\s(M(\css))|\le 2\tr D.
\end{gather}

The value on the right side of \eqref{specboun} can be computed explicitly
\begin{gather*}
\tr D:=\tr|\css_{3}|+\tr|\css_{4}'|+\tr|\css_{4}''|+\tr|\css_{4}'''|+\sum_{j=1}^{p_1-1} \tr|E_j|.
\end{gather*}
Indeed,
\begin{gather*}
|\css_{3}| =\diag\big(\big|\cb_1^2\big|, \dots, \big|\cb_{p_1}^2\big|\big), \qquad \big|\cb_j^2\big|=\diag(|b_{j,p_2}|,0,\dots,0,|b_{j,p_2}|); \\
\tr|\css_3| =\sum_{i=1}^{p_1}\tr\big|\cb_i^2\big|=2\sum_{i=1}^{p_1} |b_{i,p_2}|.
\end{gather*}
Next, since
\begin{gather*}
|G_1| =(G_1^*G_1)^{1/2}=\diag(|\a_{p_1,1}|,|\a_{p_1,2}|,\dots,|\a_{p_1,p_2}|), \\
|G_1^*| =(G_1G_1^*)^{1/2}=\diag(|\a_{p_1,p_2}|,|\a_{p_1,1}|,\dots,|\a_{p_1,p_2-1}|),
\end{gather*}
we see that
\begin{gather*} \tr|\css_4'|=2\tr|G_1|=2\sum_{i=1}^{p_2} |\a_{p_1,i}|, \qquad \tr|\css_4''|=2\tr|G_1^*|=2\sum_{i=1}^{p_2} |\a_{p_1,i}|. \end{gather*}
Clearly,
\begin{gather*} \tr|\css_4'''|=2\tr|G_2|=2\sum_{i=1}^{p_2} |a_{p_1,i}|. \end{gather*}
Finally,
\begin{gather*}
|E_j| =\diag\big(0,\dots,0,\big|\ca_j^2\big|,\big|\ca_j^2\big|,0,\dots,0\big), \\
\big|\ca_j^2\big| =\diag(|\a_{j,p_2}|,0,\dots,0,|\a_{j,p_2}|), \qquad j=1,2,\dots,p_1-1,
\end{gather*}
and so
\begin{gather*} \tr|E_j|=2\tr|\ca_j^2|=4|\a_{j,p_2}|, \qquad \sum_{j=1}^{p_1-1}\tr|E_j|=4\sum_{j=1}^{p_1-1}|\a_{j,p_2}|. \end{gather*}
A combination of the above equalities gives
\begin{gather*} \tr|D|=2\sum_{j=1}^{p_1}\bigl(|b_{j,p_2}|+2|\a_{j,p_2}|\bigr)+2\sum_{i=1}^{p_2}\bigl(|a_{p_1,i}|+2|\a_{p_1,i}|\bigr)-4|\a_{p_1,p_2}|, \end{gather*}
and hence
\begin{gather*}
|\s(J)|\le 2\tr D=4\sum_{j=1}^{p_1}\bigl(|b_{j,p_2}|+2|\a_{j,p_2}|\bigr)+4\sum_{i=1}^{p_2}\bigl(|a_{p_1,i}|+2|\a_{p_1,i}|\bigr)-8|\a_{p_1,p_2}|.
\end{gather*}

Let us note that there is nothing special in the choice of indices $p_1$, $p_2$ in the latter inequality. Indeed, it is not hard to find a unitary operator (block-shift) $W=W_{m,n}$ on $\ell^2\big(\bz,\ell^2(\bz)\big)$ so that $\hat J:=WJW^*$ is again the block-Jacobi operator~\eqref{bloper} with the shifted entries. Since $\s(J)=\s\big(\hat J\big)$, the result follows. The proof of Theorem~\ref{upbound} is complete.

Corollary~\ref{cor02} is a special case with $\a_{i,k}=0$ for $L$, and $a_{i,k}=b_{i,k}=1$ for~$H$.

\begin{Remark}\label{case2} In the case when either of $p_1$, $p_2$ equals $2$, the argument is the same with expressions \eqref{smalper}, \eqref{smalpersym} in place of \eqref{aux1}--\eqref{symb}. In particular, for $p_1=p_2=2$ we have
\begin{gather*} \frac14 |\s(J)|\le |a_{2,1}|+|a_{2,2}|+|b_{1,2}|+|b_{2,2}|+2(|\a_{1,2}|+|\a_{2,1}|+|\a_{2,2}|). \end{gather*}
\end{Remark}

\begin{Remark}\label{kruger}In \cite[Theorem~3.11]{Kr11} H.~Kr\"uger obtained a uniform upper bound for the length of the spectral bands $\ll_k$~\eqref{spect} for multidimensional periodic Schr\"odinger operators, which for $d=2$ reads
\begin{gather*} |\ll_k|\le 4\pi\left(\frac1{p_1}+\frac1{p_2}\right). \end{gather*}
The latter implies the bound for the length of the whole spectrum
\begin{gather}\label{krug}
|\s(H)|\le 4\pi(p_1+p_2).
\end{gather}
Compared to Corollary \ref{cor02}, an extra factor $\pi$ is on the right side of \eqref{krug}.
\end{Remark}

\begin{Remark}\label{gersh}
There is yet another way to obtain the upper bound for the length of $\s(J)$, based on Gershgorin's theorem \cite[Theorem 6.1.1]{hojohn}.
Indeed, let the symbol
\begin{gather*} \css(x_1,x_2)=\|s_{i,j}(x_1,x_2)\|_{i,j=1}^p. \end{gather*}
Denote by $\cg_n$ the Gershgorin disk (interval)
\begin{gather*} \cg_n(x_1,x_2):=\bigg\{\l\in\br\colon |\l-s_{nn}(x_1,x_2)|\le\sum_{j\not=n}|s_{n,j}(x_1,x_2)|, \, n=1,\dots,p\bigg\}. \end{gather*}
It is important that in our case $s_{n,n}$ and $|s_{n,j}|$ are constants, that is, do not depend on $x_1$, $x_2$, see \eqref{aux1}--\eqref{symb},
and so do the Gershgorin intervals $\cg_n$. By Gershgorin's theorem, the spectral bands~\eqref{spect}
\begin{gather*} \ll_k\subset \cg:=\bigcup_{n=1}^p \cg_n, \end{gather*}
and so
\begin{gather*} |\s(J)|\le |\cg|\le \sum_{n=1}^p |\cg_n|\le 2\sum_{i\not=j}|s_{i,j}|, \end{gather*}
the two times sum of all off-diagonal entries of the symbol $\css$. In view of \eqref{symb}, the latter leads to the bound
\begin{gather}\label{gersh1}
|\s(J)|\le \tilde R:=2\sum_{i=1}^{p_1}\sum_{j=1}^{p_2}\bigl(4|\a_{i,j}|+2|a_{i,j}|+2|b_{i,j}|\bigr).
\end{gather}

Compared to \eqref{upboun}, there is a double sum on the right side \eqref{gersh1}. It is easy to see that $\tilde R\ge\max\limits_{(m,n)\in\bt^2} R_{m,n}$,
so the bound \eqref{upboun}--\eqref{upperbound} is better than~\eqref{gersh1}.
\end{Remark}

We proceed with two examples which illustrate the optimal character of the bound in Theo\-rem~\ref{upbound}.

\begin{Example}Assume that
\begin{gather*} A_q\equiv 0, \qquad B_q=S+S^{-1}+4q I, \qquad q=1,\dots,p_1, \qquad B_{q+p_1}=B_q, \end{gather*}
$S$ is the standard shift in $\ell^2(\bz)$. In other words, we have
\begin{gather*} a_{q,r}=\a_{q,r}\equiv 0, \qquad b_{q,r}\equiv 1, \qquad c_{q,r}=4q, \qquad q=1,\dots,p_1, \qquad c_{q+p_1,r}=c_{q,r}. \end{gather*}
Now the block-Jacobi operator $J$ \eqref{genbloper} is block-diagonal, $J=\diag(B_q)_{q\in\bz}$. Since $\s(B_q)=[4q-2,4q+2]$, we see that
$\s(B_j)\cap\s(B_k)$ is at most one point set for $j,k=1,\dots,p_1$, $j\not=k$. Hence,
\begin{gather*} \s(J)=\bigcup_{q=1}^{p_1}\s(B_q)=[2,2+4p_1] \quad \Rightarrow \quad |\s(J)|=4p_1, \end{gather*}
and the factor~$4$ in Theorem~\ref{upbound} is optimal.
\end{Example}
\begin{Example}Assume that $A_q\equiv I$, $B_q\equiv B$, where $B=\diag(\b_r)_{r\in\bz}$ is a diagonal periodic matrix with
\begin{gather*} \b_r=4r, \qquad r=1,2,\dots,p_2, \qquad \b_{r+p_2}=\b_r. \end{gather*}
In this case $J=J(\{B\},\{I\})$ is unitarily equivalent to the orthogonal sum
\begin{gather*} J\simeq \bigoplus_{k=1}^{p_2}\big(J^0+4k I\big), \end{gather*}
$J^0$ is the discrete Laplacian in $\ell^2(\bz)$, so again
\begin{gather*} \s(J)=\bigcup_{k=1}^{p_2}[4k-2,4k+2]=[2,2+4p_2]\quad \Rightarrow \quad |\s(J)|=4p_2. \end{gather*}
This example again illustrates the optimal character of Theorem \ref{upbound}.
\end{Example}

\section{Examples of explicitly computed spectra}\label{s3}

The discrete version of the famous Bethe--Sommerfeld conjecture concerns the structure of the spectrum of periodic, discrete Schr\"odinger operators $H$ \eqref{disshr} (and their multidimensional analogues). It claims that for small enough potentials $\{c_{ik}\}$ such spectrum is a union of at most two closed intervals, with the gap open at the zero energy. Moreover, the spectrum is a single interval as long as at least one number $p_1$, $p_2$
is odd. The result was proved for $d=2$ in \cite{EmFi}, with a partial case for coprime periods in \cite{Kr11}, and for an arbitrary dimension $d\ge2$ in~\cite{H-J}. This result contrasts strongly with the one-dimensional case, wherein a generic $p$-periodic operator has the spectrum with $p-1$ gaps open.

It turns out that Corollary \ref{spectr} enables one to find the spectra for certain $(2,2)$-periodic discrete Schr\"odinger operators~$H$~\eqref{disshr}, with {\it not necessarily small potentials}. Indeed, by Remark~\ref{smallperiod}, the symbol is now
\begin{gather*}
\css(x_1,x_2)=
\begin{bmatrix}
\cb_1 & \tau(x_1) I_2 \\
\tau(-x_1) I_2 & \cb_2
\end{bmatrix}
\end{gather*}
with $\tau(x):=1+{\rm e}^{{\rm i}x}$,
\begin{gather*}
\cb_1=
\begin{bmatrix}
c_{11} & \tau(x_2) \\
\tau(-x_2) & c_{12}
\end{bmatrix}, \qquad
\cb_2=
\begin{bmatrix}
c_{21} & \tau(x_2) \\
\tau(-x_2) & c_{22}
\end{bmatrix}.
\end{gather*}
The characteristic polynomial of the symbol is
\begin{gather*}
\cd(\l)=\det\bigl(\css(x_1,x_2)-\l I_4\bigr)=
\begin{vmatrix}
\cb_1-\l I_2 & \tau(x_1) I_2 \\
\tau(-x_1) I_2 & \cb_2-\l I_2
\end{vmatrix}.
\end{gather*}
To compute this determinant we apply the Schur formula, which reduces determinants of order~$2n$ to ones of order~$n$ (see, e.g., \cite[Section~II.5]{Gan}). Precisely,
if $A_1$, $A_2$, $A_3$, $A_4$ are $n\times n$ matrices, and $A_1A_3=A_3A_1$, then
\begin{gather*}
\begin{vmatrix}
A_1 & A_2 \\
A_3 & A_4
\end{vmatrix} = \det(A_1A_4-A_3A_2).
\end{gather*}
Hence,
\begin{gather}\label{charpol1}
\cd(\l)=\det D(\l), \qquad D(\l)=\bigl((\cb_1-\l I_2)(\cb_2-\l I_2)-|\tau(x_1)|^2 I_2\bigr).
\end{gather}

\begin{Example}\label{ex1}
Let
\begin{gather*} c_{ij}=(-1)^{i+j}c, \qquad i,j=1,2, \qquad c>0. \end{gather*}
We find explicitly the spectrum of the corresponding operator~$H$~\eqref{disshr}:
\begin{gather*}
\s(H)=\big[{-}\sqrt{c^2+16}, -c\big]\cup \big[c,\sqrt{c^2+16}\big].
\end{gather*}
Indeed,
\begin{gather*}
(\cb_1-\l I_2)(\cb_2-\l I_2) =
\begin{bmatrix}
c-\l & \tau(x_2) \\
\tau(-x_2) & -c-\l
\end{bmatrix}
\begin{bmatrix}
-c-\l & \tau(x_2) \\
\tau(-x_2) & c-\l
\end{bmatrix} \\
\hphantom{(\cb_1-\l I_2)(\cb_2-\l I_2)}{} =
\begin{bmatrix}
\l^2-c^2+|\tau(x_2)|^2 & 2\tau(x_2)(c-\l) \\
-2\tau(-x_2)(c+\l) & \l^2-c^2+|\tau(x_2)|^2
\end{bmatrix},
\end{gather*}
and, by \eqref{charpol1},
\begin{align*}
\cd(\l) &= \big(\l^2-c^2+|\tau(x_2)|^2-|\tau(x_1)|^2\big)^2-4|\tau(x_2)|^2\big(\l^2-c^2\big) \\
&= \big(\z+|\tau(x_2)|^2-|\tau(x_1)|^2\big)^2-4|\tau(x_2)|^2\z \\
&= \big(\z-|\tau(x_2)|^2-|\tau(x_1)|^2\big)^2 -4|\tau(x_2)|^2|\tau(x_1)|^2, \qquad \z=\l^2-c^2.
\end{align*}
So, the eigenvalues of the symbol are
\begin{gather*}
\l_1(x_1,x_2) =\sqrt{c^2+(|\tau(x_2)|+|\tau(x_1)|)^2}, \\
\l_2(x_1,x_2) =\sqrt{c^2+(|\tau(x_2)|-|\tau(x_1)|)^2}, \\
\l_3(x_1,x_2) =-\l_2(x_1,x_2), \qquad \l_4(x_1,x_2)=-\l_1(x_1,x_2).
\end{gather*}
The spectral bands \eqref{spband} are
\begin{gather*} \ll_1=\big[c,\sqrt{c^2+16}\big], \qquad \ll_2=\big[c,\sqrt{c^2+4}\big], \qquad \ll_3=-\ll_2, \qquad \ll_4=-\ll_1, \end{gather*}
as claimed.
\end{Example}

\begin{Remark}The latter example is a key one in \cite{Kr11}, wherein it is shown that there is a gap in the spectrum for all $c>0$. For the similar computation of the spectrum for various lattices see \cite[Theorem~3.5]{fihan}.
\end{Remark}

\begin{Example}\label{ex2} Consider now a two-parameter family of $(2,2)$-periodic Schr\"odinger operators with
\begin{gather*}
c_{11}=-c_{22}=:c_1, \qquad c_{12}=-c_{21}=:c_2, \qquad c_1, c_2>0.
\end{gather*}
Denote $\b(x_1,x_2):=|\tau(x_2)|^2-|\tau(x_1)|^2$, so for $D$ in \eqref{charpol1} one has
\begin{gather*}
D(\l)=
\begin{bmatrix}
(\l-c_1)(\l+c_2)+\b(x_1,x_2) & -2\l\tau(x_2) \\
-2\l\tau(-x_2) & (\l-c_2)(\l+c_1)+\b(x_1,x_2)
\end{bmatrix}
\end{gather*}
and hence
\begin{gather*} \cd(\l)=\big(\l^2-c_1^2\big)\big(\l^2-c_2^2\big)+2\b(x_1,x_2)\big(\l^2-c_1c_2\big)+\b^2(x_1,x_2)-4\l^2| \tau(x_2)|^2. \end{gather*}
The roots of this biquadratic polynomial can be found explicitly. For $z=\l^2$, $\cd_1(z)=\cd(\l)$, we have
\begin{align*}
\cd_1(z) &=z^2-z\big(c_1^2+c_2^2+4|\tau(x_2)|^2-2\b(x_1,x_2)\big)+(\b(x_1,x_2)-c_1c_2)^2 \\
&= z^2-z\big(c_1^2+c_2^2+2|\tau(x_2)|^2+2|\tau(x_1)|^2\big)+(\b(x_1,x_2)-c_1c_2)^2 \\
&= z^2-2Az+B,
\end{align*}
where
\begin{gather*} A:=\frac{c_1^2+c_2^2}2+|\tau(x_1)|^2+|\tau(x_2)|^2, \qquad B:=(\b(x_1,x_2)-c_1c_2)^2\ge0. \end{gather*}
The roots $z_{\pm}$ of $\cd_1$ are
\begin{gather*}
z_{\pm}(x_1,x_2)=A\pm\sqrt{D}, \\
D= A^2-B=\left[\frac{(c_1-c_2)^2}2+2|\tau(x_2)|^2\right]\left[\frac{(c_1+c_2)^2}2+2|\tau(x_1)|^2\right]\ge0.
\end{gather*}
Hence, $0\le z_-(x_1,x_2)\le z_+(x_1,x_2)$ for each $(x_1,x_2)\in\mathbb{T}^2$, and
\begin{gather*}
\l_1(x_1,x_2) =-\l_4(x_1,x_2)=\sqrt{z_+(x_1,x_2)}, \qquad
\l_2(x_1,x_2) =-\l_3(x_1,x_2)=\sqrt{z_-(x_1,x_2)}.
\end{gather*}

We begin with the interior gap at the origin. Clearly, this gap is open if and only if
\begin{gather*} l_2=\min_{\bt^2} z_-(x_1,x_2)>0. \end{gather*}
The latter is equivalent to
\begin{gather*} \min_{\bt^2} B(x_1,x_2)=\min_{\bt^2} \bigl(2(\cos x_1-\cos x_2)-c_1c_2\bigr)^2>0. \end{gather*}
So, the interior gap at the origin is open if and only if $c_1c_2>4$.

There are another two gaps (possibly closed), symmetric with respect to the origin, which we refer to as ``exterior gaps''.
Such gaps are open if and only if
\begin{gather}\label{exam31}
\max_{\mathbb{T}^2} z_-(x_1,x_2)<\min_{\mathbb{T}^2} z_+(x_1,x_2).
\end{gather}
It is easy to see that
\begin{gather*} \min_{\mathbb{T}^2} z_+(x_1,x_2)=z_+(\pi,\pi)=\frac{c_1^2+c_2^2}2 + \frac{|c_1-c_2|}2=\max\big(c_1^2, c_2^2\big). \end{gather*}
On the other hand, it is a matter of elementary (though lengthy) calculus to check that
\begin{gather*} \max_{\mathbb{T}^2} z_-(x_1,x_2)=z_-(0,\pi)=\frac{c_1^2+c_2^2}2+4-\frac{|c_1-c_2|}2 \sqrt{(c_1+c_2)^2+16}. \end{gather*}
So, \eqref{exam31} is equivalent to
\begin{gather*}
\frac{c_1^2+c_2^2}2+4-\frac{|c_1-c_2|}2 \sqrt{(c_1+c_2)^2+16}<\max\big(c_1^2, c_2^2\big).
\end{gather*}

A simple sufficient condition for the exterior gaps to be open is
\begin{gather}\label{exam33}
c_1^2>c_2^2+8.
\end{gather}

To summarize, the conditions \eqref{exam33} and $c_1c_2>4$ ensure the existence of all three gaps, under conditions \eqref{exam33} and $c_1c_2\le4$
the exterior gaps are open. When $c_1=c_2+\ep$, $\ep>0$ small enough, and $c_2>2$ only interior gap is open.
\end{Example}

We complete with a discrete Schr\"odinger operator \eqref{disshr}, which has maximal number of gaps in its spectrum.

\begin{Example}\label{largespec}
Let $H$ be a $(p_1,p_2)$-periodic, discrete Schr\"odinger operator with the potential
\begin{gather*}
c_{i,j}=\frac{(i-1)p_2+j}{\ep} , \qquad i=1,\dots,p_1, \qquad j=1,\dots,p_2, \qquad 0<\ep<\frac18.
\end{gather*}
The symbol is now
\begin{gather*}
\css(x_1,x_2) =
\begin{bmatrix}
\cb_1 & I_{p_2} & & & {\rm e}^{{\rm i}x_1}I_{p_2} \\
I_{p_2} & \cb_2 & I_{p_2} & & \\
& \ddots & \ddots & \ddots & \\
& & I_{p_2} & \cb_{p_1-1} & I_{p_2} \\
{\rm e}^{-{\rm i}x_1}I_{p_2} & & & I_{p_2} & \cb_{p_1}
\end{bmatrix}, \\
\cb_n(x_2) =
\begin{bmatrix}
\frac{(n-1)p_2+1}{\ep} & 1 & & & {\rm e}^{{\rm i}x_2} \\
1 & \frac{(n-1)p_2+2}{\ep} & 1 & & \\
& \ddots & \ddots & \ddots & \\
& & 1 & \frac{np_2-1}{\ep} & 1 \\
{\rm e}^{-{\rm i}x_2} & & & 1 & \frac{np_2}{\ep}
\end{bmatrix}, \qquad n=1,2,\dots,p_1.
\end{gather*}
The Gershgorin intervals (see Remark~\ref{gersh}) take the form
\begin{gather*} \cg_k=\left[\frac{k}{\ep}-l, \frac{k}{\ep}+l\right], \qquad k=1,2,\dots,p, \qquad l=3 \ {\rm or} \ 4. \end{gather*}
In any event, such intervals are disjoint for $\ep<1/8$, and, by Gershgorin's theorem, each $\cg_k$ contains exactly one eigenvalue $\l_k(x_1,x_2)$ of the symbol for all $(x_1,x_2)\in\bt^2$. So, the spectral bands $\ll_k\subset\cg_k$ and, by Corollary~\ref{spectr}, there are $p-1$ gaps in the spectrum of $H$, as claimed.
\end{Example}

\subsection*{Acknowledgements}

We thank the anonymous referees for their valuable remarks which improved the paper significantly.

\pdfbookmark[1]{References}{ref}
\LastPageEnding

\end{document}